\documentclass{amsart}

\usepackage{amsthm, amssymb,latexsym,amsxtra,amscd,verbatim,ifthen}

\DeclareMathOperator{\Ker}{Ker}
\DeclareMathOperator{\Coker}{Coker}
\DeclareMathOperator{\Image}{Im}

\DeclareMathOperator{\rank}{rank} 

\DeclareMathOperator{\Supp}{Supp}
\DeclareMathOperator{\red}{red}

\DeclareMathOperator{\Spec}{Spec}

\DeclareMathOperator{\character}{char}

\DeclareMathOperator{\reg}{reg}
\DeclareMathOperator{\calHom}{\mathcal{H}\mathit{om}}
\DeclareMathOperator{\calTor}{\mathcal{T}\mathit{or}}

\newcommand{\mc}{\mathcal}

\newcommand{\ZZ}{{\mathbb Z}}

\newcommand{\NN}{{\mathbb N}}
\newcommand{\PP}{{\mathbb P}}
\newcommand{\LL}{{\mathcal L}}
\newcommand{\F}{{\mathcal F}}
\newcommand{\G}{{\mathcal G}}
\newcommand{\calE}{{\mathcal E}}

\newcommand{\calH}{{\mathcal H}}
\newcommand{\calI}{{\mathcal I}}
\newcommand{\calK}{{\mathcal K}}
\newcommand{\calN}{{\mathcal N}}

\newcommand{\calP}{{\mathcal P}}
\newcommand{\calJ}{{\mathcal J}}
\newcommand{\OO}{{\mathcal O}}

\newcommand{\gs}{{\sigma}}

\numberwithin{equation}{section}

\newtheorem{theorem}[equation]{Theorem}
\newtheorem{corollary}[equation]{Corollary}
\newtheorem{lemma}[equation]{Lemma}
\newtheorem{proposition}[equation]{Proposition}

\theoremstyle{definition}

\newtheorem{definition}[equation]{Definition}

\theoremstyle{remark}
\newtheorem{remark}[equation]{Remark}

\newtheorem{question}[equation]{Question}

\begin{document}


\title{Ample filters and Frobenius amplitude}

\author{Dennis S. Keeler}
     \thanks{ 
     Partially supported by an NSF Postdoctoral Fellowship
     and NSA Young Investigators Grant.}  
\address{Dept. of Mathematics \\ Miami University \\ Oxford, OH 45056 }
\email{keelerds@miamioh.edu}
\urladdr{http://www.users.miamioh.edu/keelerds/}
 


%
%


\begin{abstract}
Let $X$ be a projective scheme over a field.
We show that the vanishing cohomology of any sequence 
of coherent sheaves is closely related to vanishing under pullbacks
by the Frobenius morphism. We also compare various definitions
of ample locally free sheaf and show that the vanishing given by the
Frobenius morphism is, in a certain sense, the strongest possible.
Our work can be viewed as various generalizations of the Serre Vanishing
Theorem.
\end{abstract}



\maketitle

\section{Introduction}

Throughout this paper, $X$ will be a projective scheme over a field $k$.

Vanishing theorems are of great importance in algebraic geometry. They often
allow one to control the behavior of a sheaf by the behavior of its global sections.
This then allows certain algebra versus geometry connections, such
as the Serre Correspondence Theorem.

Arguably, the most important vanishing theorem is  the Serre Vanishing Theorem,
which could be stated thusly: For any coherent sheaf $\F$ on $X$, there exists
$n_0$ such that 
\begin{equation}\label{eq:vanishing}
H^q(X, \F \otimes \LL_n ) = 0
\end{equation} for $q > 0, n\geq n_0$, and
$\LL_n = \LL^{\otimes n}$ for an ample invertible sheaf $\LL$.
In noncommutative algebraic geometry, it has been extremely useful to allow
the sequence $\LL_n$ to be something other than $\LL^n$, possibly
even a sequence of non-invertible sheaves, so that one may construct
noncommutative rings from geometric data.
See \cite{AV,Keeler,Ke-filters,KRS,R,RS}.

Let $\calE$ be a fixed vector bundle and let $F$ be the absolute Frobenius
morphism on $X$ if $\character k = p > 0$.
The vanishing \eqref{eq:vanishing}
has also been studied for the sequences of locally free sheaves,
$\LL_n = S^n(\calE)$  (in which case $\calE$ is called ample) \cite{Hart-ample}
and $\LL_n = \calE^{(p^n)} = F^{*n}\calE$
(in which case $\calE$ is called cohomologically $p$-ample,
Frobenius ample, or $F$-ample) \cite{AraK,barton,gieseker}. 

Arapura recently extended this definition 
to characteristic $0$
in \cite[Definition~1.3]{AraK}, allowing him to generalize several vanishing theorems,
such as those of Kodaira-Nakano and Kawamata-Viehweg \cite[$\S$8]{AraK}. 
He also recovers
Le Potier's Vanishing Theorem by studying the ``Frobenius amplitude'' 
of a coherent sheaf. As our main theorem depends on this definition,
we give it now.
\begin{definition}\label{def:F-ample1}
Let $X$ be a projective scheme over a field $k$, and let $\F$ be a coherent sheaf.
If $\character k = p > 0$ and $F$ is the absolute Frobenius morphism, then
define $\F^{(p^n)} = F^{*n}\F$.
The \emph{Frobenius amplitude} or 
\emph{$F$-amplitude} of $\F$, $\phi(\F)$, is the smallest integer such that
for any locally free coherent sheaf $\calE$, there exists $n_0$ such that
\[
H^i(X, \calE \otimes \F^{(p^n)}) = 0,\quad i > \phi(\F), n \geq n_0.
\]
If $\character k = 0$, then $\phi(\F) \leq t$ if and only if $\phi(\F_q) \leq t$
for all closed fibers on some arithmetic thickening.
If $\phi(\F) =0$, then $\F$ is \emph{$F$-ample}.
\end{definition}

In this paper we will study the vanishing \eqref{eq:vanishing} for 
any sequence of coherent sheaves $\LL_n$. We show that the vanishing from
such a sequence is closely related to $F$-amplitude. Our main theorem is
\begin{theorem}\label{thm:smooth-case}
Let $X$ be a projective integral scheme, smooth over a perfect field $k$,
let $\{\G_n\}$ be a sequence of coherent sheaves, and
let $t \geq 0$. Then these two
statements are equivalent:
\begin{enumerate}
\item For any locally free coherent sheaf $\calE$, there exists $n_0$
such that 
\[
H^q(X, \calE \otimes \G_n)=0,\quad q > t, n\geq n_0.
\]
\item For any invertible sheaf $\calH$, there exists $n_1$
such that the $F$-amplitude $\phi(\G_n \otimes \calH) \leq t$ for $n \geq n_1$.
\end{enumerate}
\end{theorem}
(More generally, we allow the $\LL_n$ to be indexed by a filter instead
of just a sequence, because it has been useful to filter by $\NN^n$ \cite{Chan,Ke-multiple}.
See Section~\ref{S:ample filters}.)

Sections \ref{S:regularity} and \ref{S:ample filters} contain
preliminaries, though the Castelnuovo-Mumford regularity bounds
in Section~\ref{S:regularity} may be of independent interest.
The main theorem is proven in Section~\ref{S:smooth-case}.
For general projective schemes, partial generalizations are given
in Section~\ref{S:general-case}. Along the way, we obtain a
generalization of Fujita's Vanishing Theorem \cite[Theorem~1]{Fuj},
which is in turn a generalization of Serre Vanishing.
\begin{theorem}\label{thm:serre-general1}
Let $X$ be a projective scheme over a field $k$
with ample invertible sheaf $\LL$.
Let $\F$ be a coherent sheaf. Then there exists $m_0$
such that
\[
H^q(X, \F \otimes \LL^m \otimes \G) = 0,\quad q > t, m \geq m_0 
\]
for any locally free coherent $\G$ with $\phi(\G) \leq t$.
(Note that $m_0$ does not depend on $\G$.)
\end{theorem}
(See Theorems~\ref{thm:serre-smooth} and
\ref{thm:serre-general} for statements in which $\G$
is only ``nearly'' locally free.)

We then turn our attention to other ampleness properties
of a locally free sheaf.
Theorem~\ref{thm:ample-dimension-bound} shows that if $\calE$ is an ample
vector bundle (or even just ample on a complete intersection curve), 
then $\phi(\calE) < \dim X$.
Finally, in Theorem~\ref{thm:weaker}, we compare 5 possible
ways to define ampleness of a sequence of locally free sheaves.

\textbf{Acknowledgements}. We would like to thank D.~Arapura for
discussions on $F$-ampleness. We thank M.~Chardin, R.~Lazarsfeld, and J.~Sidman
for discussions on regularity. Also, H.~Brenner kindly pointed out that
the hypotheses of Theorem~\ref{thm:ample-dimension-bound} could be weakened
with no change to the proof.

\section{Subadditive regularity bounds on schemes}\label{S:regularity}

Let $X$ be a projective scheme over a field $k$ and let $\F, \G$ be coherent sheaves on $X$.
In this section we will explore how 
the Castelnuovo-Mumford regularities $\reg \F, \reg \G$, and $\reg (\F \otimes \G)$ are
related. More specifically, we generalize \cite[Proposition~1.5]{Sid},
which states that if $X = \PP^n$, then $\reg (\F \otimes \G) \leq \reg \F + \reg \G$
when $\F$ and $\G$ are ``locally free in codimension $2$.''

First, we must clarify the concept of regularity on $X$.

\begin{definition}\label{def:general-regularity}
Let $X$ be a projective scheme with an ample, globally generated invertible sheaf $\OO_X(1)$
and let $t \in \NN$.
Then a coherent sheaf $\F$ on $X$ is 
 \emph{$(m,t)$-regular} with respect to $\OO_X(1)$
if
\[
H^i(X, \F \otimes \OO_X(m-i))  = 0,\quad i > t.
\]
The minimum $m$ such that $\F$ is $(m,t)$-regular with respect to
$\OO_X(1)$ is denoted $\reg^t_{\OO_X(1)}(\F)$.
Note that if $f \colon X \to \PP^n$ is the finite morphism defined by $\OO_X(1)$,
so that $\OO_X(1) = f^*\OO_{\PP^n}(1)$, then 
$\reg^t_{\OO_X(1)}(\F) = \reg^t_{\OO_{\PP^n}(1)}(f_* \F)$.
If $\F$ is $(m,0)$-regular, we simply say $\F$ is $m$-regular.
\end{definition}

We study $(m,t)$-regularity instead of just $(m,0)$-regularity because 
it is necessary for the study of $F$-amplitude as defined in Definition~\ref{def:F-ample1}.
This concept of $(m,t)$-regularity shares a basic property with $(m,0)$-regularity.

\begin{lemma}\label{lem:greater than regularity}
Let $X$ be a projective scheme with ample, globally generated sheaf $\OO_X(1)$
and coherent sheaf $\F$. If $\F$ is $(m,t)$-regular, then $\F$ is $(n,t)$-regular
for all $n \geq m$.
\end{lemma}
\begin{proof} The proof is just as in 
\cite[p.~307, Proposition~1(i)]{Kle}.\footnote{When the base field $k$ is finite, 
one should replace
$k$ with an infinite extension. See \cite[Lemma~2.1]{Nitsure}.}
\end{proof}

In \cite[$\S$1]{AraPartial}, the related concept of the \emph{level} $\lambda(\F)$ of a 
sheaf was introduced. The level is the smallest natural number for which
$q > \lambda(\F)$ implies
\[
H^q(\F(-1)) = H^{q+1}(\F(-2)) = \dots = 0.
\]
From this we immediately have
\begin{lemma}\label{lem:partial-reg}
Let $X$ be a projective scheme with ample, globally generated sheaf $\OO_X(1)$
and coherent sheaf $\F$. Then $\lambda(\F)\leq t$ if and only if
$\F$ is $(t,t)$-regular.\qed
\end{lemma}

We can now merge, and thus slightly strengthen, \cite[Lemma~1.4]{Sid} and
\cite[Lemma~3.7]{AraK}. The reader may also confer
\cite[Section~3]{Caviglia-CM} for similar proofs for the
regularity of modules over a commutative ring.

\begin{lemma}\label{lem:complex-bounds}
Let $X$ be a projective scheme of dimension $d$
with ample, globally generated invertible sheaf $\OO_X(1)$, and let 
\[
\dots \stackrel{\phi_3}{\longrightarrow} \mc{E}_2 \stackrel{\phi_2}{\longrightarrow} 
\mc{E}_1 \stackrel{\phi_1}{\longrightarrow} 
\mc{E}_0 \longrightarrow 0
\]
be a complex of sheaves on $X$ with homology sheaves $\mc{H}_j$, for $j\geq 0$.
Let $t \in \NN$.  
Suppose that the dimension of the support of 
the higher homology of the complex satisfies
 \[\dim \Supp \mc{H}_j - j - t < 2,\quad 1 \leq j \leq d-t-2.\]  
Then \[
\reg^t_{\OO_X(1)}(\mc{H}_0) \leq \max\{ \reg^{t+j}_{\OO_X(1)}(\mc{E}_j) - j 
\colon 0 \leq j \leq d -t - 1\}.
\]
\end{lemma}
\begin{proof} 
For this proof, we assume that all regularity values are with respect to $\OO_X(1)$.
Let $m = \max\{ \reg^{t+j}(\mc{E}_j) - j \colon 0 \leq j \leq d -t - 1\}$.
Fix $i > t$. We wish to show that $H^i(X, \mc{H}_0(m-i)) = 0$.

Let $\mc{Z}_j = \Ker(\phi_j)$ and $\mc{B}_j = \Image(\phi_{j+1})$, so
that $\mc{H}_j = \mc{Z}_j/\mc{B}_j$.
The hypotheses on the $\mc{H}_j$
give us that $H^{i+j+1}(X, \mc{H}_j(m-i))=0$ for $j \geq 1$ because
$i+j+1 > \dim \Supp \mc{H}_j$. Thus 
\begin{equation}\label{eq:BjZj}
\dim H^{i+j+1}(X, \mc{B}_j(m-i)) \geq
\dim H^{i+j+1}(X, \mc{Z}_j(m-i)),\quad j \geq 1.
\footnote{The published version \cite{Ke-published} has a typo, using $\leq$ for the first $\geq$.
Thanks to S. Sierra for pointing out the error.}
\end{equation}

For $0 \leq j \leq d-t-1$, we have that $\mc{E}_j$ is $(m+j,t+j)$-regular 
by Lemma~\ref{lem:greater than regularity} 
because $m + j \geq \reg^{t+j} \mc{E}_j$. 
Thus for all $j \geq 0$, we have $H^{i+j}(X, \mc{E}_j(m-i)) = 0$
(if $j \geq d-t$, then $H^{i+j}(X, \mc{E}_j(m-i)) = 0$ because $i + j > \dim X$).
Thus
\begin{equation}\label{eq:Bj-1Zj}
\dim H^{i+j}(X, \mc{B}_{j-1}(m-i)) = \dim H^{i+j+1}(\mc{Z}_j(m-i)), \quad j \geq 1.
\end{equation}

Now $H^{i+j+1}(\mc{Z}_j(m-i))=0$ for $j \geq d -t-1$
for dimensional reasons. Then by descending induction on $j$ and 
Equations \eqref{eq:BjZj} and \eqref{eq:Bj-1Zj}, we have that 
\[ H^{i+j+1}(X, \mc{Z}_j(m-i))= H^{i+j}(X, \mc{B}_{j-1}(m-i))= 0, \quad j \geq 1.\]
We also have $H^{i}(X, \mc{Z}_0(m-i)) = 0$ because $\mc{Z}_0 = \mc{E}_0$ is $(m,t)$-regular 
by hypothesis. Thus $H^i(X, \mc{H}_0(m-i)) = 0$, as desired.
\end{proof}

The following lemma will allow us to build an appropriate complex on which to apply
Lemma~\ref{lem:complex-bounds}.

\begin{lemma}\label{lem:regularity-sequence}
Let $X$ be a projective scheme with ample, globally generated invertible sheaf $\OO_X(1)$. 
Let $\F$ be a coherent
sheaf on $X$ such that $\F$ is $m$-regular with respect to $\OO_X(1)$.
Then for any $N \geq 0$, there exist vector spaces $V_i$ and a resolution
\[
V_N \otimes \OO_X(-m-NR) \to \dots V_1 \otimes \OO_X(-m-R) \to V_0 \otimes \OO_X(-m) \to \F \to 0
\]
where $R = \max\{1, \reg^0_{\OO_X(1)}(\OO_X) \}$.
\end{lemma}
\begin{proof} The proof is exactly the same as in \cite[Corollary~3.2]{AraK},
using \cite[p.~307, Proposition~1]{Kle} for the necessary facts regarding
regularity for ample, globally generated invertible sheaves.
\end{proof}

We now generalize \cite[Proposition~1.5]{Sid}.

\begin{proposition}\label{prop:generalized-Sid}
Let $X$ be a projective scheme of dimension $d$ 
with ample, globally generated invertible sheaf $\OO_X(1)$.
Let $\F, \G$ be coherent sheaves on $X$,
and let $Y$ be the closed subscheme of $X$ where both $\F$ and $\G$ are not locally free.
If $\dim Y \leq t+ 2$, then
\[
\reg^t_{\OO_X(1)}(\F \otimes \G) \leq \reg^0_{\OO_X(1)}(\F) + \reg^t_{\OO_X(1)}(\G) + (d-t-1)(R-1)
\]
where $R = \max \{1, \reg^0_{\OO_X(1)}(\OO_X) \}$.
\end{proposition}
\begin{proof}
All regularities in this proof are with respect to $\OO_X(1)$.
Let $m = \reg^0(\F)$ and $n = \reg^t(\G)$. Form a resolution
of $\F$ as in Lemma~\ref{lem:regularity-sequence} with $N = d$. Tensoring this resolution
with $\G$, we have a complex $\mc{E}_{\bullet}$ with $\mc{E}_i = V_i \otimes \OO_X(-m-iR) \otimes \G$.
This complex has homology $\mc{H}_i = \calTor^{\OO_X}_i (\F, \G)$ and, as argued in
\cite[Proposition~1.5]{Sid}, the hypothesis on $Y$ ensures 
$\dim \Supp \mc{H}_i \leq t+2$ for $i \geq 1$.

Thus Lemma~\ref{lem:complex-bounds} applies to this complex and so
\begin{align*}
\reg^t(\mc{H}_0) &\leq \max\{ \reg^t(\mc{E}_i) - i \colon 0 \leq i \leq d - t-1\} \\
&= \max\{ m + n +i(R - 1)\colon 0 \leq i \leq d-t-1\} \\
&= m + n + (d-t-1)(R-1).
\end{align*}
Now $\mc{H}_0 = \F \otimes \G$, so we have the desired result.
\end{proof}

Note that this result is used in  \cite[Proposition~4.12]{KRS} and \cite{R,RS}.

One could make an even sharper results using $\reg^{t+i}(\G)$ for $0 \leq i \leq d-t-1$,
since $\reg^{t}(\G) \geq \reg^{t+1}(\G)$. But that becomes rather unsightly.

It seems that Proposition~\ref{prop:generalized-Sid} cannot be pushed further
to the case when $\dim Y > t+2$. Chardin has an example \cite[Example~13.6]{Chardin}
where, put in geometric terms, $X=\PP^n, \OO_X(1)$ is the degree one very ample
line bundle (so $R=1$), $\F,\G$ are ideal sheaves, and $\dim Y > 2$. In this
example, $\reg^0(\F \otimes \G) > \reg^0(\F) + \reg^0(\G)$.

\section{Ample filters of nearly locally free sheaves}\label{S:ample filters}

In this section, we take care of some technicalities regarding ample filters.
First, we must define them.

A filter $\calP$ is a   partially ordered set  such that:
\[
\text{for all } \alpha, \beta \in \calP, \text{ there exists } \gamma \in \calP
\text{ such that } \alpha < \gamma \text{ and } \beta < \gamma.
\]
Let $X$ be a projective scheme.
If a set of coherent sheaves is indexed by a filter, 
then we will call that
set a \emph{filter of sheaves}. An element of such a 
filter will be denoted $\G_\alpha$ for $\alpha \in \calP$. The indexing filter
$\calP$ will usually not be named.

\begin{definition}\label{def:ample-filter}
Let $X$ be a projective scheme over a field. Let $\calP$ be a filter.
A filter of coherent sheaves $\{ \G_\alpha \}$ on $X$, with $\alpha \in \calP$,
 will be called an 
\emph{$t$-ample filter} if for 
every locally free sheaf $\calE$, there exists
$\alpha_0$ such that
\[
H^q(X, \calE \otimes \G_\alpha) = 0, \quad q > t, \alpha \geq \alpha_0.
\]
If $\calP \cong \NN$ as filters, then a $t$-ample filter $\{ \G_\alpha \}$
is called an \emph{$t$-ample sequence}. An \emph{ample filter} is a $0$-ample filter.
\end{definition}

Note that if $X$ has characteristic $p > 0$, then the $F$-amplitude
$\phi(\G)\leq t$ if
and only if the sequence $\{ \G, \G^{(p)}, \G^{(p^2)}, \dots \}$ is a $t$-ample sequence.
Also, if $\calE$ is a semiample vector bundle, then $\calE$ is $t$-ample
if and only if $S^n(\calE)$ is a $t$-ample sequence \cite[Example~6.2.19]{PAG}.

The definition above will be useful for studying the connection between
ample filters and $F$-amplitude. However, one needs stronger vanishing theorems
to imitate the category equivalence of Serre Correspondence. (See,
for instance, \cite[Definition~5.1]{VdB-Sklyanin}.)
Fortunately, when
an ample filter consists of sheaves which are nearly locally free, we have a stronger vanishing
statement. We now make precise what we mean by nearly locally free.

\begin{definition}\label{def:lfc-n}
Let $X$ be a projective scheme, and let $\F$ be a quasi-coherent sheaf.
Then $\F$ is \emph{locally free in codimension $n$} (or \emph{lfc-$n$})
if there exists a closed subscheme $Y$ of dimension $n$ such
that $\F|_{X \setminus Y}$ is locally free. If $\F$ is locally free,
we say it is lfc-$(-1)$.
\end{definition}

It is clear that an lfc-$n$ sheaf is lfc-$(n+1)$.
Also, note that if $X$ is an integral scheme of dimension $d$, then
any coherent sheaf is lfc-$(d-1)$ \cite[Exercise~II.5.8]{Hart}.

\begin{lemma}\label{lem:lfc-2 vanishing}
Let $X$ be a projective scheme, and let $\{\G_\alpha \}$ be a $t$-ample
filter of lfc-$(t+2)$ sheaves. Then for any coherent sheaf $\F$, there
exists $\alpha_0$ such that
\[
H^q(X, \F \otimes \G_\alpha) = 0, \quad q > t, \alpha \geq \alpha_0.
\]
\end{lemma}
\begin{proof} 
Choose a very ample invertible sheaf $\OO_X(1)$.
Since $\{ \G_\alpha\}$ is a $t$-ample filter, we can find $\alpha_0$
such that $\reg^t(\G_\alpha) < -\reg^0(\F) - (\dim X - t-1)(R-1)$.
Then by Lemma~\ref{lem:greater than regularity} and
Proposition~\ref{prop:generalized-Sid}, we have
the desired vanishing.
\end{proof}

The following lemma allows us to study vanishing for lfc-$(t+1)$ sheaves.
Note that if $\F$ is lfc-$(t+1)$, then for any coherent $\G$,
$\calTor_i(\F,\G), i > 0$ is supported on a subscheme of dimension $\leq t+1$.

\begin{lemma}\label{lem:lfc-1 monomorphism}
Let $X$ be a projective scheme, let $\OO_X(1)$ be ample and generated by
global sections,
and let $0\to \calK \to  \F'' \to \F \to \F' \to 0$ be an exact
sequence of coherent sheaves. Suppose $\dim \Supp \calK \leq t+1$.
Then
\[
\reg^t(\F) \leq \max(\reg^t(\F''),\reg^t(\F')).
\]
\end{lemma}
\begin{proof}
This is immediate since $H^q(\calK(m-q)) = 0, q>t+1$ for any $m$.
\end{proof}

We will now state three standard lemmas which will help to reduce
questions about ample filters to the case where $X$ is an integral scheme.
These all follow from
applying Lemmas~\ref{lem:lfc-2 vanishing} and \ref{lem:lfc-1 monomorphism},
the methods of \cite[Exercise~III.5.7]{Hart},
and the fact that even if $\G_\alpha$ is not locally free,
the projection formula $f_*(\F \otimes f^*\G_\alpha) = f_*\F \otimes \G_\alpha$
still holds for $f$ finite \cite[Lemma~5.6]{AraK}.

\begin{lemma}\label{lem:finite morphism}
Let $f\colon Y \to X$ be a finite morphism of projective schemes and
and let $\{\G_\alpha \}$ be a filter of lfc-$(t+1)$ sheaves on $X$. If
$\{\G_\alpha\}$ is a $t$-ample filter on $X$, then
$\{ f^*\G_\alpha \}$ is a $t$-ample filter on $Y$.
If $f$ is surjective and 
$\{ f^*\G_\alpha \}$ is a $t$-ample filter on $Y$, then
$\{\G_\alpha\}$ is a $t$-ample filter on $X$.\qed
\end{lemma}

\begin{lemma}\label{lem:reduced scheme}
Let $X$ be a projective scheme with reduction $i\colon X_{\red{}} \to X$,
and let $\{ \G_\alpha \}$ be a filter of lfc-$(t+1)$ coherent sheaves on $X$.
Then $\{ \G_\alpha \}$ is a $t$-ample filter if and only if $\{ i^*\G_\alpha \}$
is a $t$-ample filter.\qed
\end{lemma}

\begin{lemma}\label{lem:irreducible components}
Let $X$ be a reduced scheme with irreducible components $f_j\colon X_j \to X, j =1 ,\dots, n$,
and let $\{ \G_\alpha \}$ be a filter of lfc-$(t+1)$ sheaves on $X$.
Then $\{ \G_\alpha \}$ is a $t$-ample filter if and only if
$\{ f^*_j\G_\alpha \}$ is a $t$-ample filter for each $j$.\qed
\end{lemma}

We will also need to change the base field so that we can work over
a perfect field.

\begin{lemma}\label{lem:base change}
Let $k \subseteq k'$ be fields, let $X$ be a projective scheme over $k$,
and let $\{ \G_\alpha \}$ be a filter of coherent sheaves on $X$.
Then $\{ \G_\alpha \}$ is a $t$-ample filter if and only if
$\{ \G_\alpha \otimes_k k' \}$ is a $t$-ample filter.
\end{lemma}
\begin{proof}
Let $f\colon X \times_k k' \to X$ be the base change, and
let $\OO_X(1)$ be a very ample line bundle on $X$. Then $f^*O_X(1)$ is very ample on $X \times k'$
\cite[II, 4.4.10]{ega-app}.
Now since $k \to k'$ is a flat morphism,
\[
 H^q(X \times_k k', f^*\G_\alpha \otimes f^*O_X(b)) = H^q(X, \G_\alpha \otimes O_X(b)) \otimes_k k'
\]
for any $b \in \ZZ$ \cite[Proposition~III.9.3]{Hart}. Then by the definition of $t$-ample 
filter and the faithful flatness
of $k \to k'$, we have the lemma. 
\end{proof}


We now wish to extend some of the basic $F$-ampleness results of \cite{AraK} 
from the locally free case to the lfc-$n$ case.

Since the definition of $F$-ampleness in characteristic zero depends on arithmetic
thickenings, we must make sure that the property lfc-$n$ behaves
well in such cases.

\begin{lemma}\label{lem:lfc-n thickening}
Let $X$ be a projective scheme over a field of characteristic $0$,
and let $\F$ be a lfc-$n$ coherent sheaf. Then
there exists an arithmetic thickening $f:\tilde{X} \to \Spec A$
such that for every fiber $\tilde{X}_x$, the sheaf $\tilde{\F}_x$ is lfc-$n$.
\end{lemma}
\begin{proof}
Let $Y$ be a closed subscheme of $X$ such that $\F|_U$ is locally free
where $U = X \setminus Y$. Then $\dim Y \leq n$ by hypothesis.
We can find a thickening $\tilde{Y} \stackrel{i}{\to} \tilde{X} \stackrel{f}{\to} \Spec A$
such that $i$ is a closed immersion and each fiber of $f \circ i$ has dimension $\leq n$.
We can shrink the thickening so that $\F$ has a thickening $\tilde{\F}$ and
$\tilde{\F}|_{\tilde{U}}$ is locally free, where $\tilde{U} = \tilde{X} \setminus \tilde{Y}$.
Then for every $y \in \Spec A$, we have $\tilde{\F}|_{\tilde{U}_y}$ locally free
and $\dim \tilde{Y}_y \leq n$. Thus $\tilde{\F}_y$ is lfc-$n$, as desired.
\end{proof}

Because the study of $F$-ampleness relies on pulling back by the Frobenius
morphism $F$, we also need to see how the lfc-$n$ property behaves
under pullback.

\begin{lemma}\label{lem:lfc-n pullback}
Let $X$ be a projective scheme, and let $\F$ be a lfc-$n$ coherent sheaf on $X$.
If $f\colon X' \to X$ is a quasi-finite morphism (that is,
a morphism with finite fibers), then $f^* \F$ is lfc-$n$.
\end{lemma}
\begin{proof}
Let $Y \subset X$ be a closed subscheme such that $\dim Y \leq n$ and
$\F_{X \setminus Y}$ is locally free. Then
 $f^*(\F_{X \setminus Y}) = f^*\F_{X' \setminus f^{-1}(Y)}$ is locally free
 and $\dim f^{-1}(Y) \leq n$ since $f$ is quasi-finite.
\end{proof}

We now turn our attention to $F$-amplitude (see Definition~\ref{def:F-ample1}).
Note that the definition only requires vanishing when tensoring with
a locally free sheaf. However, if $\F$ is locally free in codimension $2$,
we can get more, generalizing \cite[Lemma~2.2]{AraK}.

\begin{lemma}\label{lem:F-ample gives vanishing}
Let $X$ be a projective scheme over a field of characteristic $p > 0$,
and let $\F$ be a lfc-$(t+2)$ coherent sheaf. Then for any coherent sheaf $\G$,
there exists $n_0$ such that
\[
H^i(X, \G \otimes \F^{(p^n)}) = 0,\quad i > \phi(\F), n \geq n_0.
\]
\end{lemma}
\begin{proof}
The proof uses the same methods as Lemma~\ref{lem:lfc-2 vanishing}.
\end{proof}

\section{Ample filters and $F$-ampleness: smooth case}\label{S:smooth-case}

Our main goal in this section is to prove Theorem~\ref{thm:smooth-case}.
For this section, $X$ will be a projective variety (that is, an integral scheme), 
smooth over a perfect
field $k$.

We will first need a lemma which is a slight generalization of
\cite[Theorem~5.4]{AraPartial}. The concept of ``sufficiently ample line bundle''
is defined in \cite[$\S$1]{AraPartial}. For our purposes, we need
only know that given any ample line bundle $\LL$, there exists $n_0$
such that $\LL^n$ is sufficiently ample for $n \geq n_0$ 
\cite[Theorem~1.2]{AraPartial}.

\begin{lemma}\label{lem:t-ample-gives-F-ample}
Let $X$ be a projective variety over a perfect field of characteristic $p>0$.
Let $\OO_X(1)$ be a sufficiently ample line bundle, and 
let $\calE, \F$ be coherent sheaves, with at least
one of $\calE, \F$ locally free.
If $q > \lambda(\calE(-\dim X))$ and if $p^N \geq \reg^0(\F)$, then
\[
H^q(X, \calE^{(p^N)} \otimes \F) = 0.
\]
\end{lemma}
\begin{proof}
The proof proceeds exactly as in \cite[Theorem~5.4]{AraPartial},
where it is assumed both $\calE, \F$ are locally free.
We need only make sure that  a certain exact sequence
remains exact when tensoring with
$\calE \otimes \F$.

More specifically, Arapura sets $f'=F^n$, the $n$th power of the absolute
Frobenius. Then $f':X \to X$ is not $k$-linear, but since 
$k$ is perfect, we can let $X' = X \times_{\Spec k} \Spec k$
over the $p^n$th power map of $k$. There is then
a $k$-linear map $f$ and natural map $g$:
\[
X \stackrel{f}{\to} X' \stackrel{g}{\to} X,
\]
where $f' = g\circ f$. Set $\calE' = g^*\calE$.

We then let $\gamma: X \to X'\times X$ be the morphism for which $p_1\circ \gamma=f$
and $p_2 \circ \gamma = id$.
Arapura constructs a locally free resolution of $\gamma_*\OO_X$, which
he tensors with $\calE' \boxtimes \F = p_1^*\calE' \otimes p_2^*\F$.
If we can verify that $\calTor_i(\calE' \boxtimes \F, \gamma_*\OO_X)=0, i>0$,
then we are done. Since one of $\calE', \F$ is locally free, it
suffices to verify that $\calTor_i(p_1^*\calE', \gamma_*\OO_X)=0$
and $\calTor_i(p_2^*\F, \gamma_*\OO_X)=0$.

Let $0 \to \calK \to \G \to \calE \to 0$ be exact with $\G$ locally free.
Set $\calK' = g^*\calK, \G' = g^*\G$. Then since $p_1$ is flat, there
is a short exact sequence
\[
0 \to p_1^* \calK' \to p_1^*\G' \to p_1^*\calE' \to 0
\]
and hence an exact sequence
\[
0 \to \calTor_1(p_1^*\calE', \gamma_*\OO_X) \to 
p_1^* \calK' \otimes  \gamma_*\OO_X \to p_1^*\G' \otimes  \gamma_*\OO_X
\to p_1^*\calE'\otimes  \gamma_*\OO_X \to 0.
\]

However, since $X,X'$ are
 smooth (and hence regular) \cite[Prop.~VII.6.3]{AltKle}
 and $f$ is finite, 
we have that $f$ is flat
\cite[Corollary~V.3.6]{AltKle}. Since $\gamma$ is also finite, there
is an exact sequence
\[
0 \to 
\gamma_*(f^*\calK') \to \gamma_*(f^*\G')
\to \gamma_*(f^*\calE') \to 0. 
\]
Since $\gamma$ is a finite morphism, the projection formula 
holds for any coherent sheaf \cite[Lemma~5.6]{AraK}. That is,
$p_1^*\calE'\otimes  \gamma_*\OO_X \cong \gamma_*(\gamma^*p_1^*\calE')
\cong  \gamma_*(f^*\calE')$. So we have
\[
0 \to \calTor_1(p_1^*\calE', \gamma_*\OO_X) \to 
\gamma_*(f^*\calK') \to \gamma_*(f^*\G')
\to \gamma_*(f^*\calE') \to 0.
\]
Since all these exact sequences were derived by application
of functors, we have $\calTor_1(p_1^*\calE', \gamma_*\OO_X) = 0$.

We also have $\calTor_i(p_1^* \calK',  \gamma_*\OO_X) \cong 
\calTor_{i+1}(p_1^* \calE',  \gamma_*\OO_X)$ for $i>0$. Allowing $\calK'$
to play the role of $\calE'$ and using induction, we have
$\calTor_{i}(p_1^* \calE',  \gamma_*\OO_X) = 0$ for $i > 0$.

The case of $\calTor_{i}(p_2^* \F,  \gamma_*\OO_X) = 0$ is
similar.
\end{proof}

We immediately have
a generalization of \cite[Corollary~5.5]{AraPartial}.
\begin{corollary}\label{cor:F-amplitude-less-than-partial}
Let $X$ be a smooth projective variety defined over a perfect
field $k$ of arbitrary characteristic. Let $\OO_X(1)$
be sufficiently ample and let $\calE$ be a coherent 
sheaf on $X$. Then the Frobenius amplitude $\phi(\calE)$
satisfies 
\[
\phi(\calE) \leq \lambda(\calE(-\dim X)).
\]
\end{corollary}
\begin{proof} Using an arithmetic thickening, we may assume $k$
has characteristic $p > 0$. Then letting $\F$ be an arbitrary
locally free sheaf, the result follows from Lemma~\ref{lem:t-ample-gives-F-ample}.
\end{proof}

We can now prove one direction of Theorem~\ref{thm:smooth-case}.
\begin{theorem}\label{thm:smooth-case1}
Let $X$ be a projective variety, smooth over a perfect field $k$. 
Let $\{\G_\alpha\}$ be a $t$-ample filter. Then for
any locally free coherent sheaf $\calE$, there exists $\alpha_0$
such that $\phi(\G_\alpha \otimes \calE) \leq t$ for $\alpha \geq \alpha_0$.
\end{theorem}
\begin{proof} Since $\{\G_\alpha\}$ is a $t$-ample filter, there
exists an index $\alpha_0$ such that $\lambda(\G_\alpha \otimes \calE(-\dim X)) \leq t$
for $\alpha \geq \alpha_0$. We then apply Corollary~\ref{cor:F-amplitude-less-than-partial}.
\end{proof}

We now turn to proving the other half of Theorem~\ref{thm:smooth-case}.
The main tool is the following exact sequence.
When $k$ is algebraically closed, the surjectivity of the trace map
below is mentioned, but not proven, in \cite[$\S$1]{Fuj}. It is presumably well-known,
but we provide a proof for lack of reference.

\begin{lemma}\label{lem:flat-sequence}
Let $X$ be a projective variety, smooth over a perfect field of 
characteristic $p > 0$. Let $F:X \to X$ be the absolute Frobenius morphism.
Then there is an exact sequence of coherent, locally free sheaves
\[
0 \to \calK \to F_*\omega \to \omega \to 0
\]
where $\omega$ is the canonical sheaf of $X$.
\end{lemma}
\begin{proof}
Since $X$ is smooth
  (and hence regular) \cite[Proposition~VII.6.3]{AltKle}
 and $F$ is finite, 
we have that $F$ is finite and flat
\cite[Corollary~V.3.6]{AltKle}.
Thus $F_* \OO_X$ is a locally free sheaf.
Consider $\phi:\OO_X \to F_*\OO_X$
with the morphism $\phi$ locally given by $a \mapsto a^p$. The morphism $\phi$ 
remains injective
locally at $x \in X$
upon tensoring with any residue field $k(x)$, so the cokernel of $\phi$ is also
locally free \cite[Exercise~II.5.8]{Hart}.
Let $\omega = \wedge^{\dim X} \Omega^1_{X/k}$
\cite[Theorem~VII.5.1]{AltKle}.
Dualizing via \cite[Exercises~III.6.10, 7.2]{Hart}, we have the desired
short exact sequence of locally free sheaves.
\end{proof}

We can now exploit the method of \cite[Lemma~5.8]{Fuj2}, \cite[$\S$1]{Fuj}
to prove a version of Fujita's Vanishing Theorem \cite[Theorem 1]{Fuj}
which is a generalization of Serre Vanishing.

\begin{theorem}\label{thm:serre-smooth}
Let $X$ be a projective variety, smooth over a perfect field $k$
with ample invertible sheaf $\LL$.
Let $\calE$ be a locally free coherent sheaf. Then there exists $m_0$
such that
\[
H^q(X, \calE \otimes \LL^m \otimes \G) = 0, \quad q > t, m \geq m_0 
\]
for any coherent $\G$ with $\phi(\G) \leq t$.
(Note that $m_0$ does not depend on $\G$.)
\end{theorem}
\begin{proof}
First suppose that the characteristic of $k$ is $p > 0$. Let
$0 \to \calK \to F_* \omega \to \omega \to 0$ be the exact sequence
of Lemma~\ref{lem:flat-sequence}. Let $q > t$ be given. 
Suppose by descending induction on $q$ that
\[
H^{q+1} (X,  \calK \otimes \LL^m \otimes \G) = 0,\quad m \geq m_1, \phi(\G) \leq t.
\]
(This is certainly true for $q \geq \dim X$ because
$\dim X$ is the cohomological dimension of $X$.)

Let $h^q(\F) = \dim_k H^q(X,\F)$.
Then we have
\[
h^q(F_*\omega \otimes \LL^m \otimes \G)
= h^q(F_*(\omega \otimes \LL^{pm} \otimes \G^{(p)}))
\geq h^q(\omega \otimes \LL^m \otimes \G)
\]
for $m \geq m_1$ and \emph{any} $\G$ with $\phi(\G) \leq t$. 
But if $\phi(\G) \leq t$, then $\phi(\LL^m \otimes \G^{(p^n)}) \leq t$ for any 
$n \geq 0, m\geq 0$ \cite[Theorem~4.5]{AraK}.
So then $h^q(\omega \otimes \LL^{p^em} \otimes \G^{(p^e)})
\geq h^q(\omega \otimes \LL^m \otimes \G)$ for all $e \geq 0$.
Since $\phi(\LL^m \otimes \G^{(p^n)}) \leq t$, these dimensions equal $0$
for $e$ sufficiently large. But 
then $ h^q(\omega \otimes \LL^m \otimes \G) = 0, q > t, m \geq m_1$.

Now suppose the characteristic of $k$ is $0$. Then a version of the Kodaira
Vanishing Theorem holds. That is, if $\phi(\G) \leq t$, then
$H^q(X, \omega \otimes \G) = 0, q > t$ \cite[Corollary~8.6]{AraK}.

Now let the characteristic of $k$ be arbitrary. Let $\calE$ be a locally
free coherent sheaf. Let $\OO_X(1)$ be very ample. By taking $m_2$
sufficiently large, we can make $\reg^0(\calE \otimes \omega^{-1} \otimes \LL^{m_2})<0$.
And by the above arguments, by taking $m_3$ sufficiently large,
$\reg^t(\omega \otimes \LL^{m} \otimes \G) < -(\dim X - t - 1)(R-1)$
where $R = \max\{1, \reg^0(\OO_X)\}$ and $m \geq m_3$.
Then by Proposition~\ref{prop:generalized-Sid}, we have
$\reg^t(\calE \otimes \LL^{m} \otimes \G) < 0$ for $m \geq m_2 + m_3$.
This immediately gives the theorem.
\end{proof}

We can now prove the other half of Theorem~\ref{thm:smooth-case}.

\begin{theorem}\label{thm:smooth-case2}
Let $X$ be a projective variety, smooth over a perfect field $k$.
Let $\{ \G_\alpha \}$ be a filter of coherent sheaves. 
Suppose that for any locally free sheaf $\calE$, there exists
$\alpha_0$ such that the Frobenius amplitude $\phi(\G_\alpha \otimes \calE) \leq t$
for $\alpha \geq \alpha_0$. Then $\{\G_\alpha \}$ is a $t$-ample sequence.
\end{theorem}
\begin{proof}
Let $\LL$ be an ample invertible sheaf and 
let $\F$ be a coherent locally free sheaf. Then by Theorem~\ref{thm:serre-smooth},
there exists $m$ such that $H^q(X, \F \otimes \LL^m \otimes \G)=0$
for all coherent $\G$ with $\phi(\G) \leq t$. For sufficiently large
$\alpha$, we have $\phi(\LL^{-m} \otimes \G_\alpha) \leq t$, and we are done.
\end{proof}

\section{Ample filters and $F$-ampleness: general case}\label{S:general-case}

In this section, we will obtain some generalizations of Theorem~\ref{thm:smooth-case}
to the case of $X$ a projective scheme over a general field.
Unfortunately, we are unable to obtain Theorem~\ref{thm:smooth-case}
in full generality. If $X$ is not smooth, then the Frobenius morphism
is not flat \cite[Theorem~2.1]{Kunz}. We also lose the invertibility
of $\omega_X$ \cite[Theorem~VII.5.3]{AltKle}.
These facts lead to weaker results.

The following is a partial generalization of Theorem~\ref{thm:smooth-case1}.
Note that the proof only works for $0$-ample filters.

\begin{proposition}\label{prop:ample filter implies F-ample}
Let $X$ be a projective scheme over a field $k$, and let $\{\G_\alpha \}$ be an ample
filter of lfc-$2$ sheaves. Then for any invertible sheaf $\calH$,
 there exists $\alpha_0$ such that
$\calH \otimes \G_\alpha$ is $F$-ample for all $\alpha \geq \alpha_0$.
\end{proposition}
\begin{proof} 
Since $\{ \calH \otimes \G_\alpha \}$ is an ample filter, we may assume that
$\calH = \OO_X$.
Let $d = \dim X$,
let $\OO_X(1)$ be a very ample invertible sheaf, and let $R = \max\{1, \reg(\OO_X)\}$.
 Since $\{ \G_\alpha \}$ is an ample filter, there exists
$\alpha_0$ such that $\reg(\G_\alpha) < -(d-1)R$ for $\alpha \geq \alpha_0$.
Then by Lemma~\ref{lem:regularity-sequence}, there is a locally free resolution
$\calE_\bullet \to \G_\alpha \to 0$ where $\calE_i$ is a direct sum of very ample
invertible sheaves, $i = 0,\dots,d-1$.

If $\character k = p > 0$, then pulling back $\calE_\bullet \to \G_\alpha \to 0$
via $F^{*n}$ for large $n$, and using Lemma~\ref{lem:complex-bounds}, we see that
$\lim_{n \to \infty} \reg(F^{*n} \G_\alpha) = -\infty$.
Thus Proposition~\ref{prop:generalized-Sid} shows that for any coherent locally free sheaf
$\calE$, the higher cohomology
of $\calE \otimes F^{*n}\G_\alpha$ vanishes for large $n$, so $\G_\alpha$ is $F$-ample.
If $\character k = 0$, then one can choose an arithmetic thickening so
that $\tilde{\calE}_\bullet \to \tilde{\G}_\alpha \to 0$ is still a resolution
by direct sums of very ample invertible sheaves, at least for $i = 0, \dots, d-1$.
Then the above argument on each closed fiber of the thickening shows that $\G_\alpha$
is $F$-ample for $\alpha \geq \alpha_0$.
\end{proof}

For the remainder of the section, we will prove a partial generalization
of Theorem~\ref{thm:smooth-case2}. The main idea of the proof is to
reduce to the smooth case via alteration of singularities. Thus,
we must study the higher direct images of a projective morphism $f:X \to Y$.

\begin{lemma}\label{lem:higher-direct}
Let $f:X \to Y$ be a projective morphism of projective schemes over a field.
Let $\LL$ be an $f$-ample invertible sheaf on $X$ and let $\F$ be a coherent
sheaf on $X$. Then there exists $m_0$ such that
\begin{align}\label{eq:higher-direct}
f_*(\F \otimes \LL^m) \otimes \G &\cong f_*(\F \otimes \LL^m \otimes f^*\G), \notag\\
R^if_*(\F \otimes \LL^m) \otimes \G &= R^if_*(\F \otimes \LL^m \otimes f^*\G)=0,
\end{align}
for all $i>0, m \geq m_0$, and coherent $\G$ on $Y$. (Note that $m_0$ is 
independent of $\G$.)

Thus also
\[
H^q(Y, f_*(\F \otimes \LL^m) \otimes \G) \cong H^q(X, \F \otimes \LL^m \otimes f^*\G)
\]
for $q \geq 0, m \geq m_0$, and all coherent $\G$ on $Y$.
\end{lemma}
\begin{proof}
Our proof is similar to the proof of \cite[Lemma~5.7]{AraK}.
We proceed by descending induction on $i$. By Grothendieck's
Vanishing Theorem \cite[Theorem~III.2.7]{Hart},
the higher direct images are identically $0$ for large $i$.
So assume we have \eqref{eq:higher-direct} for $i+1$.

Let $\G$ be a coherent sheaf on $Y$
and let $\calE_1 \to \calE_0 \to \G \to 0$ be a resolution of $\G$
by locally free sheaves. Let $\F(m) = \F\otimes \LL^m$.
Then for some $m_0$, we have
a commutative diagram
\begin{equation*}
\minCDarrowwidth20pt
\begin{CD}
R^if_*(\F(m)) \otimes \calE_1 
@>>> R^if_*(\F(m)) \otimes \calE_0
@>>> R^if_*(\F(m)) \otimes \G \to 0  \\
@VVV @VVV @VVV \\
R^if_*(\F(m) \otimes f^*\calE_1) 
@>>> R^if_*(\F(m) \otimes f^*\calE_0)
@>>> R^if_*(\F(m) \otimes f^*\G) \to 0
\end{CD}\end{equation*}
for $m \geq m_0$.
The first two vertical arrows are isomorphisms by the usual
projection formula \cite[Exercise~III.8.3]{Hart}
and thus the third one is as well. 

Now if $i > 0$, then we may increase $m_0$ so that
$R^if_*(\F \otimes \LL^m) = 0$ for $m \geq m_0$ \cite[Theorem~1.7.6]{PAG}.
This completes the proof of \eqref{eq:higher-direct}.
The last statement of the lemma follows immediately from 
\cite[Exercise~III.8.1]{Hart}.
\end{proof}

Since the ampleness of a line bundle is not preserved under pull-backs,
$F$-ampleness cannot be preserved either. However, if $\LL$ is ample,
then $f^*\LL$ is nef, and hence $\calH \otimes f^*\LL$ is ample
for any ample line bundle $\calH$. This does generalize.

\begin{lemma}\label{lem:F-ample-pullbacks}
Let $X,Y$ be projective schemes over a perfect field $k$ and 
let $f:X \to Y$ be a projective morphism.
Let $\F$ be a lfc-$(t+1)$ coherent sheaf on $Y$ with $\phi(\F) \leq t$.
Then for any ample line bundle $\calH$ on $X$, we have
$\phi(\calH \otimes f^*\F) \leq t$.
\end{lemma}
\begin{proof} 
Suppose $\character k = p>0$.
Since $k$ is perfect, the absolute Frobenius morphism
is finite. Thus $\phi(\calH \otimes f^*\F) = \phi(\calH^{p^n} \otimes
f^*\F^{(p^n)})$ by Lemma~\ref{lem:finite morphism}. 
Thus we may assume that $\calH$ is very ample.
Define $\OO(1) = \calH$.

Let $\calE$ be a locally free coherent sheaf on $X$. 
Any ample line bundle is $f$-ample 
\cite[Proposition~1.7.10]{PAG}. So
by Lemma~\ref{lem:higher-direct}, there exists $m_0$ such that
\[
H^q(Y, f_*(\calE(m-q)) \otimes \G) \cong H^q(X, \calE(m-q) \otimes f^*\G)
\]
for all $q\geq 0, m \geq m_0$ and coherent $\G$ on $Y$.
Since $\phi(\F) \leq t$ and $\F$ is lfc-$(t+2)$, there exists $n_0$ such that
\[
H^q(Y, f_*(\calE(m_0-q)) \otimes \F^{(p^n)}) = 0
\]
for $q > t, n \geq n_0$ by Lemma~\ref{lem:F-ample gives vanishing}.
Thus 
$\calE \otimes f^*\F^{(p^n)}$ is $(m,t)$-regular for $m \geq m_0, n\geq n_0$.
For $n_1 \geq n_0$ sufficiently large, we have $p^{n_1} \geq m_0+\dim X$.
Thus for $n \geq n_1, q > t$, we have
\[
H^q(X, \calE(p^n)\otimes f^*\F^{(p^n)}) = 0.
\]
Thus $\phi(\calH \otimes f^*\F) \leq t$.

If $\character k = 0$, the result follows immediately.
\end{proof}

We may now prove another generalization of Serre-Fujita Vanishing
by reducing to Theorem~\ref{thm:serre-smooth}.

\begin{theorem}\label{thm:serre-general}
Let $X$ be a projective scheme over a field $k$
with ample invertible sheaf $\LL$.
Let $\calE$ be a coherent sheaf. Then there exists $m_0$
such that
\[
H^q(X, \calE \otimes \LL^m \otimes \G) = 0,\quad  q > t, m \geq m_0
\]
for any coherent lfc-$(t+1)$ $\G$ with $\phi(\G) \leq t$.
(Note that $m_0$ does not depend on $\G$.)
\end{theorem}
\begin{proof}
Let $\Sigma = \{ \LL^m \otimes \G : \phi(\G) \leq t, \G \text{ lfc}-(t+1) \}$.
Partially order $\Sigma$ by $\LL^{m_1} \otimes \G_1 < \LL^{m_2} \otimes \G_2$
if $m_1 < m_2$. If we can show that $\Sigma$ is an ample filter, 
then we are done by Lemma~\ref{lem:lfc-2 vanishing}.
By Lemmas~\ref{lem:reduced scheme}, \ref{lem:irreducible components},
and \ref{lem:base change}, we may assume that $X$ is an integral scheme
over a perfect field.

Let $f: \tilde{X} \to X$ be an alteration of singularities,
so that $\tilde{X}$ is smooth over $k$ and $f$ is surjective
\cite{deJong-Alterations}. Let $\calH$ be an ample line bundle
on $\tilde{X}$. Then by Theorem~\ref{thm:serre-smooth}
and Lemma~\ref{lem:F-ample-pullbacks}, there
exists $m_0$ such that
\[
H^q(\tilde{X}, \calH^m \otimes f^*\G) = 0
\]
for $q > t, m \geq m_0, \phi(\G) \leq t$, with $\G$ lfc-$(t+1)$.
Then by Lemma~\ref{lem:higher-direct}, for some $m_1 \geq m_0$, we have
\[
H^q(X, f_*(\calH^m) \otimes \G) = 0
\]
for $q > t, m \geq m_1, \phi(\G) \leq t$, with $\G$ lfc-$(t+1)$.

Fix $m$. Since $f$ is surjective,
we have $\Supp f_*(\calH^m) = X$. Thus the sheaf $\calHom(f_*(\calH^m),\OO_X)$ 
contains a non-zero (and hence surjective) local homomorphism at the generic point.
For $N$ sufficiently large,
\[
\calHom(f_*(\calH^m),\OO_X) \otimes \LL^N \cong 
\calHom(f_*(\calH^m),\LL^N)
\]
is generated by global sections.
Therefore there exists a global homomorphism $\psi:f_*(\calH^m) \to \LL^N$
which is surjective at the generic point. Thus $\Coker \psi$ is a 
torsion sheaf.

Now descending induction on $q$ and noetherian induction
applied to $\Coker \psi$ gives $H^q( \LL^n \otimes \G) = 0$
for $q > t, n \geq n_0, \phi(\G) \leq t,$ with $\G$ lfc-$(t+1)$.
If $\calE$ is any coherent sheaf on $X$, we may find
a resolution of $\calE$ by direct sums of $\LL^{-j}$. The
theorem then follows by Lemma~\ref{lem:complex-bounds}.
\end{proof}

We now immediately find the following partial generalization of
Theorem~\ref{thm:smooth-case2}, using the same method of proof.

\begin{corollary}\label{cor:general-case2}
Let $X$ be a projective scheme over a field $k$.
Let $\{ \G_\alpha \}$ be a filter of lfc-$(t+1)$ coherent sheaves. 
Suppose that for any invertible sheaf $\calH$, there exists
$\alpha_0$ such that the Frobenius amplitude $\phi(\G_\alpha \otimes \calH) \leq t$
for $\alpha \geq \alpha_0$. Then $\{\G_\alpha \}$ is a $t$-ample filter.\qed
\end{corollary}

\section{$F$-amplitude of an ample vector bundle}\label{S:ample vector}

Let $\calE$ be an ample vector bundle. If $\character k = 0$, then
$\phi(\calE) < \rank(\calE)$ \cite[Corollary~6.6]{AraK}. We will now derive
another bound on $\phi(\calE)$, which is independent of the characteristic of $k$.
First, we need a lemma.

\begin{lemma}\label{lem:very-ample-bounds}
Let $X$ be a projective scheme, let $H$ be a very ample Cartier divisor,
and let $\calE$ be a vector bundle. Then
\[
\phi(\calE_H) \leq \phi(\calE) \leq \phi(\calE_H) + 1.
\]
\end{lemma}
\begin{proof}
The first inequality follows from Lemma~\ref{lem:finite morphism}.

If $\character k = p > 0$, 
for any $b \in \ZZ$ there exists $n_0$ such that $\OO_X(b) \otimes \calE_H^{(p^n)}$
is $(0, \phi(\calE_H))$-regular for $n \geq n_0$. Thus there are exact sequences
\begin{multline*}
0 = H^q(O_H(b+m-q) \otimes \calE_H^{(p^n)}) 
\to H^{q+1}(O_X(b+m-q-1) \otimes \calE^{(p^n)}) 
\\ \to H^{q+1}(O_X(b+m-q) \otimes \calE^{(p^n)})
\to H^{q+1}(O_H(b+m-q) \otimes \calE_H^{(p^n)}) = 0
\end{multline*}
for $q > \phi(\calE_H), m\geq 0, n \geq n_0$. By Serre Vanishing, 
$H^{q+1}(O_X(b+m-q) \otimes \calE^{(p^n)}) = 0$ for $m \gg 0$. 
So by descending induction on $m$ and \cite[Corollary~2.3]{AraK},
we have
$\phi(\calE) \leq \phi(\calE_H) + 1$.
The case of $\character k = 0$ is then immediate.
\end{proof}

It is now an easy matter to obtain our bound on $\phi(\calE)$ for ample $\calE$.
This generalizes 
\cite[Proposition~5.4]{AraK}. 

\begin{theorem}\label{thm:ample-dimension-bound}
Let $X$ be a projective scheme over a field $k$ with
 $\dim X > 0$, and let $\calE$ be a vector bundle
 which is ample on some complete intersection curve.
Then $\phi(\calE) < \dim X$. In particular, this is
true for ample $\calE$.
\end{theorem}
\begin{proof}
We may assume that $k$ is algebraically closed \eqref{lem:base change} and that
$X$ is reduced \eqref{lem:reduced scheme}, irreducible  
\eqref{lem:irreducible components}, and normal
\eqref{lem:finite morphism}.
If $\dim X = 1$, then the claim is \cite[Proposition~5.4]{AraK}.
If $\dim X > 1$, then induction on $\dim X$ and Lemma~\ref{lem:very-ample-bounds}
yields the result.
\end{proof}

We can expect no better result. On $\PP^n$, the tangent bundle
$\mathcal{T}_{\PP^n}$ is ample (and even $p$-ample, see Definition~\ref{def:p-ample ample}),
yet $\phi(\mathcal{T}_{\PP^n}) = n-1$ by \cite[Theorem~5.5]{AraK}, because
$H^{n-1}(\mathcal{T}(-n-1)) \neq 0$ \cite[Example~5.9]{AraK}.

\section{Other ampleness conditions}\label{S:other ampleness}

In this section we compare various alternative definitions of ``ample filter'' for locally free sheaves
and indicate how they relate to the Frobenius morphism. We remind the reader of the definition of $p$-ample 
and ample locally free sheaf.

\begin{definition}\label{def:p-ample ample}
Let $X$ be a projective scheme over a field $k$, and let $\calE$ be a locally free coherent sheaf.
If $\character k = p > 0$, then $\calE$ is \emph{$p$-ample} if for any coherent
sheaf $\F$, there exists $n_0$ such that $\F \otimes \calE^{(p^n)}$ is generated by
global sections for $n \geq n_0$. If $\character k = 0$, we say $\calE$ is $p$-ample
if $\calE$ is $p$-ample on every closed fiber of some arithmetic thickening.

If the sequence $\{ S^n(\calE) : n \in \NN \}$ is an ample sequence, then $\calE$
is an \emph{ample} locally free sheaf.
\end{definition}

We now state some progressively weaker ``ampleness'' properties on filters of locally free sheaves.

\begin{theorem}\label{thm:weaker}
Let $X$ be a projective scheme, and let $\{ \calE_\alpha \}$ be a filter of
locally free coherent sheaves. Consider the following properties.
\begin{enumerate}
\item\label{F-ample1} For any coherent sheaf $\F$, there exists $\alpha_0$ such that
$H^q(\F \otimes \calE_\alpha)=0$ for $q>0, \alpha \geq \alpha_0$.
\item\label{F-ample2} For any invertible sheaf $\calH$, there exists $\alpha_0$ such
that $\calH \otimes \calE_\alpha$ is an $F$-ample locally free sheaf for $\alpha \geq \alpha_0$.
\item\label{p-ample1} For any coherent sheaf $\F$, there exists $\alpha_0$ such that
$\F \otimes \calE_\alpha$ is generated by global sections for $\alpha \geq \alpha_0$.
\item\label{p-ample2} For any invertible sheaf $\calH$, there exists $\alpha_0$ such
that $\calH \otimes \calE_\alpha$ is a $p$-ample locally free sheaf for $\alpha \geq \alpha_0$.
\item\label{ample} For any invertible sheaf $\calH$, there exists $\alpha_0$ such
that $\calH \otimes \calE_\alpha$ is an ample locally free sheaf for $\alpha \geq \alpha_0$.
\end{enumerate}
Then (\ref{F-ample1}) $\Longleftrightarrow$ (\ref{F-ample2}) $\implies$ (\ref{p-ample1}) $\implies$
(\ref{p-ample2}) $\implies$ (\ref{ample}).
\end{theorem}
\begin{proof} \eqref{F-ample1} $\Leftrightarrow$ \eqref{F-ample2} is 
Proposition~\ref{prop:ample filter implies F-ample} and 
Corollary~\ref{cor:general-case2}.
For \eqref{F-ample1} implies \eqref{p-ample1}, we may choose $\alpha_0$ such that $\F \otimes \calE_\alpha$
is $0$-regular, and hence generated by global sections, for $\alpha \geq \alpha_0$.

For \eqref{p-ample1} implies \eqref{p-ample2},
let $\LL$ be an ample invertible sheaf, and let $\calH$ be an arbitrary invertible sheaf. Choose
$\alpha_0$ such that $\LL^{-1} \otimes \calH \otimes \calE_\alpha$ is generated by global
sections for $\alpha \geq \alpha_0$. Then $\calH \otimes \calE_\alpha$ is the quotient of
a finite direct sum of copies of $\LL$. 
It is easy to see in characteristic $p$
 that $\oplus \LL$ is $p$-ample and that quotients of $p$-ample sheaves are $p$-ample
\cite[Proposition~6.4]{Hart-ample}.
 If $\character k = 0$, the surjection 
$\oplus \LL \twoheadrightarrow \calH \otimes \calE_\alpha$ can be extended to some arithmetic
thickening. Hence the closed fibers of
$\tilde{\calH} \otimes \tilde{\calE}_\alpha$
are $p$-ample, so $\calH \otimes \calE_\alpha$ is $p$-ample for $\alpha \geq \alpha_0$.

For \eqref{p-ample2} implies \eqref{ample}, it is known in characteristic $p > 0$ that
$p$-ample bundles are ample \cite[Proposition~6.3]{Hart-ample}. In characteristic $0$, each of the closed fibers of
a thickening $\tilde{\calH} \otimes \tilde{\calE}_\alpha$ is ample, hence 
$\calH \otimes \calE$ is ample by \cite[$\mathrm{III}_1$, 4.7.1]{ega-app}.
\end{proof}

In certain special cases, all these statements are equivalent.

\begin{proposition}\label{prop:reverse}
Let $X$ be a projective scheme, and let $\{ \calE_\alpha \}$ be a
filter of locally free coherent sheaves. Suppose one of the following
holds.
\begin{enumerate}
\item\label{curve} $\dim X \leq 1$,
\item\label{invertible} Each $\calE_\alpha$ is an invertible sheaf,
\item\label{symmetric} The filter is  $\{ S^n(\calE) : n \in \NN \}$
for some fixed locally free coherent sheaf $\calE$.
\end{enumerate}
Then all the statements of Theorem~\ref{thm:weaker} are equivalent.
\end{proposition}
\begin{proof} If $\dim X \leq 1$ or $\calE_\alpha$ is invertible, then
ample implies $F$-ample by
Theorem~\ref{thm:ample-dimension-bound} or \cite[Theorem~1.3]{Ke-filters}.
Therefore, \ref{thm:weaker}\eqref{ample} implies \ref{thm:weaker}\eqref{F-ample2}.
If $\{ S^n(\calE) : n \in \NN \}$ satisfies \ref{thm:weaker}\eqref{ample},
then $S^m(\calE)$ is ample for some $m$. So then $\calE$ is ample 
\cite[Proposition~2.4]{Hart-ample}. But the vanishing in \ref{thm:weaker}\eqref{F-ample1}
can be taken as a definition of an ample locally free sheaf, at least when $X$ is proper
\cite[Proposition~3.3]{Hart-ample}.
\end{proof}

We note that the equivalence of statements still holds when the $\calE_\alpha$
are invertible sheaves and $X$ is proper over a commutative noetherian ring
\cite[Theorem~1.3]{Ke-filters}.

In general, the implications are not reversible, 
except possibly \eqref{p-ample1} $\implies$ \eqref{p-ample2}.

\begin{remark}\label{prop:no reverse}
Let $X =\PP^2$ over an algebraically closed field of characteristic $p > 0$.
Then there exists a locally free coherent sheaf $\calE$ such that
$\{ \calE^{(p^n)} : n \in \NN \}$ satisfies \ref{thm:weaker}(\ref{p-ample1}), but
not \ref{thm:weaker}(\ref{F-ample2}).
There also exists a locally free coherent sheaf $\F$ such that
$\{ \F^{(p^n)} : n \in \NN \}$ satisfies \ref{thm:weaker}(\ref{ample}), but
not \ref{thm:weaker}(\ref{p-ample2}).
\end{remark}
\begin{proof}
Let $\calE$ be $p$-ample, but not $F$-ample. Then the first claim is satisfied by
definition. Let $\F$ be ample, but not $p$-ample. If
the sequence $\{ \F^{(p^n)} \}$ satisfies \ref{thm:weaker}\eqref{p-ample2}, then
$\F^{(p^n)}$ is $p$-ample for some $m$. But then $\F$ is $p$-ample 
\cite[Proposition~6.4]{Hart-ample}, a contradiction. However, 
$\{ \F^{(p^n)} \}$ satisfies \ref{thm:weaker}\eqref{ample} \cite[Proposition~3.1]{barton}.
Such $\calE, \F$ do exist on $\PP^2$ \cite{gieseker}. An important specific example
is the tangent bundle $\mathcal{T}_{\PP^n}$ for $n \geq 2$. This
bundle is $p$-ample, but not $F$-ample \cite[Example~5.9]{AraK}.
\end{proof}

\begin{question} Does \ref{thm:weaker}(\ref{p-ample2}) imply \ref{thm:weaker}(\ref{p-ample1})
in general? The answer is affirmative when $X$ is smooth \cite{Ke-FujitaConj}.
\end{question}

\bibliographystyle{amsplain}

\begin{thebibliography}{EGA}

\bibitem[AK]{AltKle}
A.~Altman and S.~Kleiman, \emph{Introduction to {G}rothendieck duality theory},
  Springer-Verlag, Berlin, 1970.

\bibitem[A1]{AraK}
D. Arapura, 
\emph{Frobenius amplitude and strong vanishing theorems for vector bundles}, 
with an appendix by Dennis S. Keeler, Duke Math. J., 
\textbf{121} (2004), no.~2, 231--267.

\bibitem[A2]{AraPartial}
\bysame, \emph{Partial Regularity and Amplitude}, Amer. J. Math. \textbf{128}
(2006), no.~4, 1025--1056.
  
\bibitem[AV]{AV}
M.~Artin and M.~Van~den Bergh, \emph{Twisted homogeneous coordinate 
rings}, J. Algebra \textbf{133} (1990), no.~2, 249--271.
  
\bibitem[B]{barton}
C.~M. Barton, \emph{Tensor products of ample vector bundles in
  characteristic $p$}, Amer. J. Math. \textbf{93} (1971), 429--438. 
  
\bibitem[Ca]{Caviglia-CM}
G. Caviglia, \emph{{Bounds on the {C}astelnuovo-{M}umford regularity of tensor
  products}}, Proc. Amer. Math. Soc. \textbf{135} (2007), 1949--1957.
  
\bibitem[Chn]{Chan}
D. Chan, \emph{Twisted multi-homogeneous coordinate rings}, J. Algebra
  \textbf{223} (2000), no.~2, 438--456. 
  
\bibitem[Chr]{Chardin}
M. Chardin, \emph{Some results and questions on {C}astelnuovo-{M}umford
  regularity}, Syzygies and Hilbert functions, Lect. Notes Pure Appl. Math.,
  vol. 254, Chapman \& Hall/CRC, Boca Raton, FL, 2007, pp.~1--40.

\bibitem[D]{deJong-Alterations}
A.~J. de~Jong, \emph{Smoothness, semi-stability and alterations}, Inst. Hautes
  \'Etudes Sci. Publ. Math. (1996), no.~83, 51--93. 

\bibitem[F1]{Fuj2}
T. Fujita, \emph{Semipositive line bundles}, {J. Fac. Sci. Univ. Tokyo Sect.
  IA Math.} \textbf{30} (1983), no.~2, 353--378. 
  
\bibitem[F2]{Fuj}
\bysame, \emph{Vanishing theorems for semipositive line bundles},
  Algebraic geometry (Tokyo-Kyoto, 1982), Springer, Berlin, 1983, 
pp.~519--528.

\bibitem[G]{gieseker}
D. Gieseker, \emph{$p$-ample bundles and their {C}hern classes}, Nagoya
  Math. J. \textbf{43} (1971), 91--116. 

\bibitem[EGA]{ega-app} 
A.~Grothendieck, \emph{\'{E}l\'ements de g\'eom\'etrie alg\'ebrique}, Inst.
  Hautes \'Etudes Sci. Publ. Math. (1961, 1966), no.~8, 11, 28.

\bibitem[H1]{Hart-ample} 
R. Hartshorne, \emph{Ample vector bundles}, Inst. Hautes \'Etudes Sci. Publ.
  Math. (1966), no.~29, 63--94.

\bibitem[H2]{Hart}
\bysame, \emph{Algebraic geometry}, Graduate Texts in Math., no.~52,
  Springer-Verlag, New York, 1977.
  
\bibitem[Ke0]{Ke-published}
D.~S. Keeler, \emph{Ample filters and Frobenius amplitude},
 J. Algebra \textbf{323} (2010), no.~ 10, 3039–-3053.

\bibitem[Ke1]{Keeler}
D.~S. Keeler, \emph{Criteria for $\sigma$-ampleness}, J. Amer. Math. Soc.
  \textbf{13} (2000), no.~3, 517--532.
  
\bibitem[Ke2]{Ke-multiple}  
\bysame, \emph{Noncommutative ampleness for multiple divisors}, J.
  Algebra \textbf{265} (2003), no.~1, 299--311. 
  
\bibitem[Ke3]{Ke-filters}
\bysame, \emph{Ample filters of invertible sheaves}, J. Algebra
  \textbf{259} (2003), no.~1, 243--283.
  
\bibitem[Ke4]{Ke-FujitaConj}
\bysame, \emph{{F}ujita's {C}onjecture and {F}robenius amplitude}, 
Amer. J. Math. \textbf{130} (2008), no.~5, 1327--1336.

\bibitem[KRS]{KRS}
D. S. Keeler, D. Rogalski, J. T. Stafford, \emph{Na\"ive noncommutative blowing up}, Duke Math J.,
\textbf{126} (2005), no.~3, 491--546. 

\bibitem[Kle]{Kle}
S.~L. Kleiman, \emph{Toward a numerical theory of ampleness}, Ann. of Math.
  (2) \textbf{84} (1966), 293--344.

\bibitem[Ku]{Kunz} E. Kunz,
\emph{Characterizations of regular local rings for characteristic $p$},
Amer. J. Math. \textbf{91} (1969) 772--784.

\bibitem[L]{PAG}
R. Lazarsfeld, \emph{Positivity in algebraic geometry. {I}, {II}}, Springer-Verlag, Berlin,
  2004.


\bibitem[N]{Nitsure} N. Nitsure,
\emph{Construction of Hilbert and Quot schemes}. Fundamental algebraic geometry, 105–137,
Math. Surveys Monogr., 123, Amer. Math. Soc., Providence, RI, 2005. 

\bibitem[R]{R} D. Rogalski, \emph{GK-dimension of birationally commutative surfaces},
Trans. Amer. Math. Soc. \textbf{361}  (2009), no. 11, 5921--5945.

\bibitem[RS]{RS} D. Rogalski, J. T. Stafford, 
\emph{Na\"ive noncommutative blowups at zero-dimensional schemes},
J.\ Algebra  \textbf{318}  (2007), no. 2, 794--833.

\bibitem[S]{Sid}
J. Sidman, \emph{On the {C}astelnuovo-{M}umford regularity of products of ideal sheaves}, 
Adv. Geom., \textbf{2} (2002), no.~3, 219--229.

\bibitem[VdB]{VdB-Sklyanin}
M.~Van~den Bergh, \emph{A translation principle for the four-dimensional
  {S}klyanin algebras}, J.\ Algebra \textbf{184} (1996), 435--490.

\end{thebibliography}
\providecommand{\bysame}{\leavevmode\hbox to3em{\hrulefill}\thinspace}

\end{document}